\documentclass[graybox]{svmult}

\usepackage{booktabs}     

\usepackage{type1cm}        
\usepackage{makeidx}         
\usepackage{graphicx}        
\usepackage{multicol}        
\usepackage[bottom]{footmisc}

\usepackage{newtxtext}       %
\usepackage[varvw]{newtxmath}       

\makeindex             

\begin{document}

\title*{Investment Decisions for Perfect and Imperfect Competition in Ireland's Electricity Market}
\titlerunning{Electricity Market Investment in Ireland: Perfect and Imperfect Competition}
\author{Davoud Hosseinnezhad \orcidID{0000-0002-4011-776X} \\ Mel T. Devine \orcidID{0009-0009-2393-3026}
\\ Seán McGarraghy \orcidID{0000-0001-9906-8635}}
\institute{Davoud Hosseinnezhad \at University College Dublin, Dublin, Ireland, \email{Dhoseinnejad@yahoo.com}
\and Mel T. Devine \at College of Business, University College Dublin, Dublin, Ireland \email{Mel.Devine@ucd.ie}
\and Seán McGarraghy \at College of Business, University College Dublin, Dublin, Ireland \email{Seán.McGarraghy@ucd.ie}}
%
%
\maketitle

\vspace{-6em}
\abstract{This paper employs a game-theoretic approach to analyze investment decisions in Ireland’s electricity market. It compares optimal electricity investment strategies among energy generators under a perfect competition framework with an imperfect Nash-Cournot competition. The model incorporates market price based on competition among generators while accounting for supply capacity of each firm and each technology, along with the System Non-Synchronous Penetration (SNSP) constraint to reflect operational limitations in the renewable energy contribution to the power system. Both models are formulated as single-objective function optimization problems. 
Furthermore, unit commitment constraints are introduced to the perfect competition model allowing the model to incorporate binary decision variables to capture energy unit scheduling decisions of online status, startup, and shutdown costs. The proposed models are evaluated under three different demand test cases, using Ireland's electricity generation projections for 2023-2033. The results highlight key differences in investment decisions, $CO_2$ emissions, and the contribution of renewable technologies in perfect and imperfect competition structures. The findings provide managerial insights for policy makers and stakeholders, supporting optimal investment decisions and generation capacity planning to achieve Ireland’s long-term energy objectives.}

\textbf{Keywords} Energy Market, Competition Models, Unit Commitment Constraints
\vspace{-1em}

\section{Introduction}
\label{sec:1}
\vspace{-1em}
Game theory is widely used in the Operations Research literature on energy markets to analyze strategic interactions among market participants \cite{pozo2017bilevel}. A primary distinction in this domain lies between \textbf{perfect competition} and \textbf{imperfect competition}, which influence market dynamics, pricing mechanisms, and participant behavior. These classifications determine whether energy suppliers act as price-takers or price-makers, shaping their strategic decision-making and overall market outcomes \cite{tomic2024perfect}.

In a perfectly competitive energy market, all participants are price-takers, meaning they accept the market price as given and have no ability to influence it. This model assumes a large number of suppliers offering homogeneous products, leading to a perfectly elastic demand curve for each supplier \cite{rozas2024comparative}. Perfect competition models are often used to analyze markets in which individual energy suppliers have small market shares and prices are determined solely by aggregate supply and demand \cite{budzinski2024perfect, devine2024analysing}. Perfect competition models can be formulated as \textbf{cost minimization problems} or equivalently \textbf{social welfare maximization problems} \cite{egging2020solving}. 

Unlike models based on Nash-Cournot competition or equilibrium concepts, \textbf{cost minimization} does not incorporate strategic competition among firms, as it assumes a centralized approach where demand is determined externally \cite{devine2012rolling}. This method has been widely applied in various market models \cite{sundarajoo2023under, faranda2007load}. \textbf{Social welfare maximization} is also applied to energy market problems to address perfectly competitive market structures \cite{egging2020solving, devine2016rolling, ramnath2021social}.

Incorporating integer variables into energy market models allows the representation of key operational constraints, such as unit commitment and generation scheduling \cite{kercci2020impact}. The \textbf{perfect competition with unit commitment problem} is a critical aspect of power system dispatch, determining which generators should be online at any given time while accounting for operational constraints such as startup and online costs \cite{devine2023cournot}. Previous studies, such as \cite{tuohy2009unit}, have explored the integration of renewable generation into unit commitment models, while more recent works (e.g \cite{wu2023novel, ali2024multi, zhang2023two}) focus on developing efficient solution techniques for these models. See \cite{montero2022review} for a comprehensive review of unit commitment problems and approaches.

Imperfect competition, in contrast, captures market structures where certain participants can influence prices, acting as price-makers. These suppliers strategically adjust production quantities or prices to maximize profits while anticipating competitors’ responses. Imperfect competition is common in oligopolistic markets, where some companies have important market power allowing them to influence market prices through their decisions and are aware of this influence \cite{pozo2017bilevel}. Imperfect competition models are also used to model capacity restrictions and transmission limits \cite{ruiz2011equilibria, gonzalez2022transmission}. These models pose greater challenges in incorporating integer variables. However, some exceptions exist, such as the work by \cite{devine2023cournot}. Imperfect competition can be categorized into different types based on how participants make decisions. \textbf{Cournot competition}, \textbf{Bertrand competition}, and \textbf{Supply Function Equilibria (SFE)} are three main imperfect competition models discussed in the literature \cite{pozo2017bilevel}.

\textbf{Cournot competition models} formulate firms competing by deciding on energy generation quantities. It is foundational in studying oligopolistic markets, where firms compete based on output rather than price. Each firm selects its generation level to maximize its profit, considering how its decisions influence market prices and the strategies of competitors. Equilibrium is reached when no firm can unilaterally increase its profit by adjusting its production. Unlike perfectly competitive markets, where firms have no market power and take the market-clearing price as given, Cournot competition assumes that firms recognize their influence on prices. Market prices are determined by the inverse demand function, reflecting the aggregate supply decisions of all firms \cite{trebicka2014imperfect, hobbs2002linear}. Nash-Cournot games can be solved as a quadratic programming problem, assuming the principles of symmetry are met \cite{egging2020solving}.

While Cournot competition models are preferred when adjusting production quantities is complex, \textbf{Bertrand models} suit markets with flexible pricing \cite{gabriel2012complementarity}. \cite{saukh2020modelling, cartuyvels2024market} applied the Bertrand model in electricity markets.
\textbf{SFE models} examine electricity markets where firms submit supply curves based on price and quantity. SFE outcomes resemble Cournot equilibria during peak demand and Bertrand equilibria in off-peak periods \cite{pozo2017bilevel}. See \cite{wei2023supply, mokhtari2021impact} for a review of SFE models for electricity markets.


Both perfect competition and imperfect competition models can be formulated and solved using either \textbf{single-level optimization approaches} or \textbf{bi-level optimization approaches}. In bi-level optimization, decisions are made in a hierarchical sequence of leader-follower decisions. This framework, known as \textbf{Stackelberg equilibrium}, is widely applied in power system models, such as strategic bidding \cite{lei2023data} and generation investment \cite{weibelzahl2020optimal}. A more complex extension, the multiple-leader-multiple-follower equilibrium, involves interactions between multiple decision-makers at both levels, often modeled using \textbf{Equilibrium Problem with Equilibrium Constraints} (EPEC) \cite{luna2023regularized,devine2023strategic}. However, EPEC models are computationally challenging, especially for large-scale systems.

The contribution of this work can be expressed by reviewing relevant studies. \cite{devine2023cournot} consider perfect and imperfect competition, and also unit commitment model, using detailed data on the Irish electricity system. However, their work does not incorporate long-term investment decisions. \cite{helisto2019backbone} include unit commitment constraints and investment decisions with detailed Irish system data but do not model imperfect competition. \cite{devine2024analysing} incorporate investment decisions in both perfect competition and imperfect frameworks, but their analysis is limited to toy data and does not compare outcomes with a unit commitment model.

Building on these studies, we employ a single-level optimization approach and develop game-theoretic models for Ireland’s electricity market to analyze decisions for both perfect and imperfect competition models. Furthermore, the perfect competition model is extended by incorporating unit commitment constraints, where startup and no-load costs are modeled using binary variables. The key contribution of this study lies in integrating long-term investment decisions within both perfect and imperfect competition models, using detailed Irish market data and benchmarking the outcomes against a unit commitment model. Such an approach has not been seen in the literature to date.

The remainder of this paper is structured as follows. Section \ref{sec:2} introduces the Irish energy market problem, detailing generation units and technologies. Section \ref{sec:3} presents the mathematical models for the proposed single objective maximization problem with conjectural variation parameter for perfect and imperfect competition problem, followed by the extended perfect competition model with unit commitment constraints. Section \ref{sec:4} validates these models under Ireland’s projected 2026 energy demand test cases followed by the results and conclusions in Section \ref{sec:5}.
\vspace{-1em}

\section{Proposed Model}
\label{sec:2}
\vspace{-1em}
This paper examines optimal investment and electricity generation decisions in Ireland's electricity market, where multiple electricity-generating firms ($f = 1,2, \dots, F$) compete by determining the quantity of energy they supply to meet market demand. Each firm operates a set of current generation units, indexed by $u'_f = 1,2, \dots, u'_F$, where $u'_F$ represents the number of existing units owned by firm $f$. The firms use various energy generation technologies, and each unit employs a distinct technology. For instance, a firm may operate one unit generating electricity from gas and another from coal. These units can be in different geographical locations but remain under the ownership and capacity portfolio of the same firm.

While firms operate their existing units, they also have the option to invest in new generation units. However, investment is restricted to specific energy sources, namely gas, coal, hydro, oil, wind, and solar that is indexed by $u' = 1,2, \dots, 6$, respectively. Consequently, the total set of units for each firm consists of both existing and new units, expressed as $u_f = u'_f + u'$, where $u'_f$ represents the set of current units, and $u'$ represents the new generation units to invest. 

Each firm makes two primary decisions at each time $t$ ($t=1,2, \dots, T$) under three capacity factors scenarios $s$ ($s = 1,2,3$) which affect the maximum availability of solar and wind technologies throughout the year. The probability of each scenario occurring is formulated as $P_s$, where the sum of all scenario probabilities satisfies $\sum_{s=1}^{3} P_s = 1$.
The decision variables (all non-negative) of the proposed model are: (i) The generation decision, denoted as $q_{f,u_f,t,s}$, which represents the electricity generated by the unit $u_f$ of firm $f$ at time $t$ under scenario $s$. (ii) The investment decision, denoted as $\text{Inv}_{f,u_f}$, which determines the firm’s investment in each unit. Since investment in existing units is not allowed, we enforce $\text{Inv}_{f,u_f} = 0$ for all existing units ($u'_f$).

\begin{figure}[!hbt]
    \centering
    \includegraphics[width=10 cm]{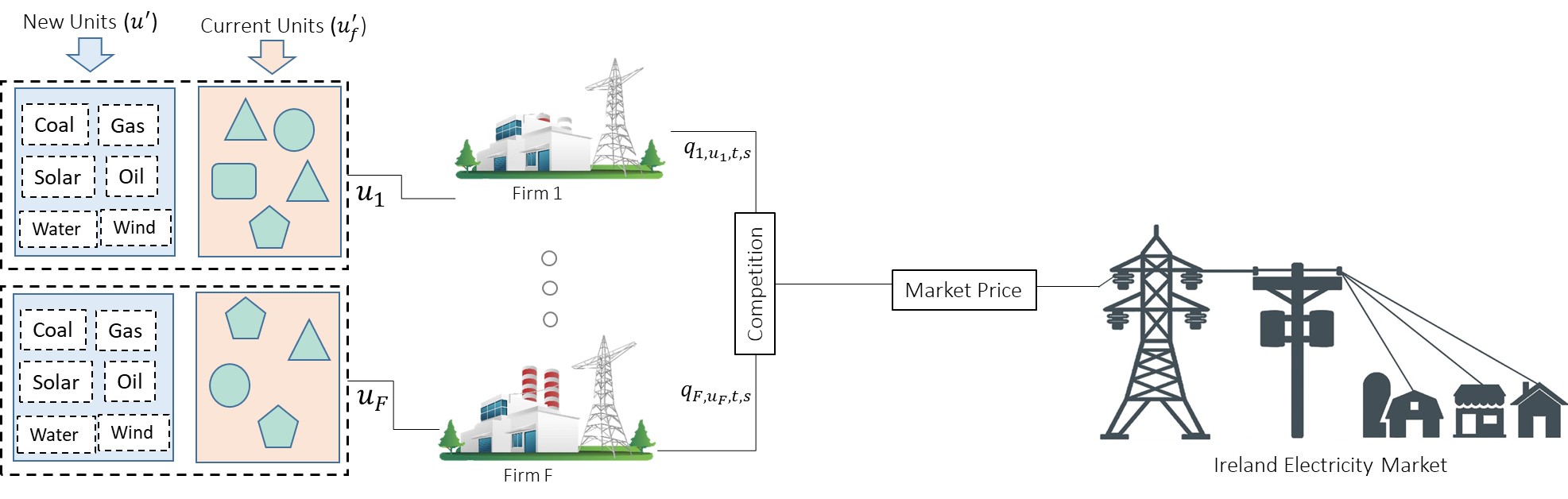}
    \caption{The Proposed Electricity Generation Network in Ireland}
    \label{fig1}
\end{figure}

A schematic representation of the proposed model is illustrated in Figure \ref{fig1}. The market price of the electricity at each time and under each scenario (\( \pi_{t,s} \)) is determined by the inverse demand function, which in turn influences firms' decisions. The market price function at time \( t \) and scenario $s$ is given by \( \pi_{t,s} = A_t - B \sum_{f,u_f} q_{f,u_f,t,s} \), where \( A_t \) represents the price when total supply is zero, and \( B \) is a constant representing the slope of the demand curve. 
\vspace{-1em}
\vspace{-1em}

\section{Mathematical Formulation}
\label{sec:3}
\vspace{-1em}
Each firm's objective is to maximize its profit, which consists of the revenue from selling electricity in the market, minus the production costs, and minus investment costs for establishing additional generation units. The optimization problem for each firm $f$ is formulated as follows:
\vspace{-0.5em}
\begin{align}
\label{equation1}
\max_{q,{Inv}} & \sum_{t,s} P_s \times W_t \times \sum_{u_f} \left((\pi_{t,s}- C_{f,u_f}) \cdot q_{f,u_f,t,s} -C_{u_f}^{\text{Inv}} \cdot \text{Inv}_{f,u_f}\right) \\
\text{subject to:} \nonumber \\
& q_{f,u_f,t,s} \leq CF_{f,u_f,t,s}\cdot (Q_{f,u_f}^{\max} + \text{Inv}_{f,u_f})  \quad \forall u_f,t,s , \label{equation2} 
\end{align}

in which $W_t$ is the weighting factor of the profit function at each time $t$. $C_{f,u_f}$ represents the marginal cost of firm $f$ producing energy from a unit $u_f$. $C_{u_f}^{\text{Inv}}$ is the marginal investment cost for establishing a new generation unit of type $u_f$.

Each unit \( u_f \) operated by firm \( f \) has a maximum generation capacity denoted as \( Q_{f,u_f}^{\max} \). This represents the upper limit on the energy output that an existing unit can produce at the time \(t\). However, for new units to be established (\( u' \)), we assume that \( Q_{f,u'}^{\max} = 0 \). The capacity constraint in Equation (\ref{equation2}) ensures that the quantity of energy produced by firm $f$ from unit $u_f$ at time $t$ does not exceed the total available capacity of that unit at that time. The available capacity in this equation is scaled by the capacity factor \( CF_{f,u_f,t,s} \), reflecting that renewable technologies in units may not operate at full capacity. \( CF_{f,u_f,t,s} \) is between 0 and 1. $q_{f,u_f,t,s}$ and $\text{Inv}_{f,u_f}$ are non-negative variables.

To solve the proposed mathematical model, we can employ the Karush-Kuhn-Tucker conditions to formulate the equilibrium condition between the firms. 
However, this paper uses the approach in \cite{egging2020solving}, which involves reformulating the model as a single optimization model and can be solved as a quadratically constrained problem. This approach's advantage lies in its ability to directly incorporate market power exertion within the objective function and reduce computation time. 
\vspace{-1.5em}

\subsection{Recasting the Competition as a Single Optimization Problem}
\label{subsec:3-1}
\vspace{-1em}
In this section, we frame competition between generating firms as a single-maximization problem that does not require the deviation or implementation of the Karush-Kuhn-Tucker conditions. Accordingly, the objective function is the sum of the profit functions of each firm, augmented by an appropriate conjectural variation term that accounts for the market power of the firms \cite{egging2020solving}. The objective function of the proposed single optimization problem is expressed in Equation (\ref{equation9}).
\vspace{-0.5em}
\begin{align}
\label{equation9}
\text{Max} = & \sum_{t,s} P_s \times W_t \times  \Bigg[ 
    \sum_{f,u_f} \Big( 
          (\pi_{t,s}- C_{f,u_f}) \cdot q_{f,u_f,t,s} -C_{u_f}^{\text{Inv}} \cdot \text{Inv}_{f,u_f}
    \Big) \nonumber \\
    & + \frac{1}{2} \cdot B \cdot \Big( \sum_{f,u_f} q_{f,u_f,t,s} \Big)^2
    - \theta \cdot \Big[\frac{1}{2} \cdot B \cdot \sum_{f} \Big( \sum_{u_f} q_{f,u_f,t,s} \Big)^2 
\Big] \Bigg] 
\end{align}
Subject to the Equation (\ref{equation2}) for all firms and the following constraint:
\begin{align}
\sum_{f,u_f^{\text{renewable}}} q_{f,u_f,t,s} \leq 0.75 \sum_{f,u_f}q_{f,u_f,t,s} \quad \forall t,s \label{equation11}
\end{align}

The parameter \(\theta\) in the objective function captures the degree of market power exerted by the firms \cite{egging2020solving}. We continue to analyze optimal decisions of the firms for two extreme values of \(\theta\) (\(\theta = 0\) which represents perfect competition, and \(\theta = 1\) which represents Nash-Cournot imperfect competition). Equation \ref{equation11} is the System Non-Synchronous Penetration (SNSP) constraint that represents a physical limit of the power network, restricting to 75\% the maximum share of non-synchronous generation that can be accommodated at any moment while ensuring system stability.
\vspace{-1em}

\subsection{Perfect competition with unit commitment constraints}
\label{subsec:3-2}
\vspace{-1em}
In this extended formulation, we extend the single optimization model of perfect competition ($\theta=0$) by introducing binary variables that represent the status of each generation unit. Accordingly, the binary variable $\text{on}_{f, u_f, t, s}$ and $\text{sd}_{f, u_f, t, s}$ indicate whether each unit is online or shutdown at time $t$, while the binary variable $\text{su}_{f, u_f, t, s}$ states if a unit remains online at time $t$. The modified objective function now includes associated startup costs ($C^{su}_{f, u_f}$) and online operation ($C^{su}_{f, u_f}$) in the objective function. The modified perfect competition model with unit commitment constraints, also named social welfare optimization with unit commitments (SW-UC), is as follows.
\vspace{-0.5em}
\begin{align}
\label{equation14}
\text{Max} & \text{ SW-UC} =  \sum_{t,s} P_s \times W_t \times  \Bigg[ 
    \sum_{f,u_f} \Big( 
        (\pi_{t,s}- C_{f,u_f}) \cdot q_{f,u_f,t,s} -C_{u_f}^{\text{Inv}} \cdot \text{Inv}_{f,u_f}
    \Big) + \nonumber \\
&\frac{1}{2} B \cdot \Big( \sum_{f,u_f} q_{f,u_f,t,s} \Big)^2
    \Bigg] - \sum_{f,u_f} \Big( C^{on}_{f, u_f} \cdot \text{on}_{f,u_f,t,s} + C^{su}_{f, u_f} \cdot \text{su}_{f,u,t,s} \Big) 
\end{align}

Subject to Equation (\ref{equation11}) and constraints (\ref{equation15}) to (\ref{equation18}). 

\begin{align}
q_{f,u_f,t,s} \leq \text{CF}_{f,u_f,t,s} \cdot \left(\text{on}_{f,u_f,t,s} \cdot Q_{f,u_f}^{\max} + \text{Inv}_{f,u_f}\right) \quad \forall f, u_f, t,s \label{equation15} \\
q_{f,u_f,t,s} \geq \text{on}_{f,u_f,t,s} \cdot Q_{f,u_f}^{\min} \quad \forall f, u_f, t,s \label{equation16} \\
\text{su}_{f,u_f,t,s} - \text{sd}_{f,u_f,t,s} = \text{on}_{f,u_f,t,s} - \text{on}_{f,u_f,t-1,s} \quad \forall f, u_f, t, s \label{equation17}\\
\text{su}_{f,u_f,t,s} + \text{sd}_{f,u_f,t,s} \leq 1 \quad \forall f, u_f, t, s \label{equation18}
\end{align}
Equation (\ref{equation16}) states the minimum generation level for each unit if it is online. Equations (\ref{equation17} and \ref{equation18}) enforce the logic between the startup, shutdown, and online status of units. Accordingly, a unit can only start if it was previously off and is now being brought online. Also, a unit cannot start up and shut down at the same time.
\vspace{-1em}

\section{Input Data}
\label{sec:4}
\vspace{-1em}
The proposed models for the Irish electricity market include all 16 generation firms in Ireland, each owning a number of existing units. 
However, investment is restricted to specific energy sources, namely
gas, coal, hydro, oil, wind, and solar. The analysis is carried out on an hourly basis, with each time period \( t \) representing one hour, over a span of three weeks. Each week corresponds to a week in winter, spring/autumn, and summer seasons. The objective function for this span is multiplied by an appropriate importance factor $W_t$ so that the decision variables represent annual optimal investment and generation decisions. Each of the proposed models has been run for three demand test cases (Low, Median, and High) available on Ireland's \textit{Ten-Year Generation Capacity Statement} \cite{eirgrid2023generation}. The demand intercept values follow from \cite{devine2023cournot} but are calibrated to these three test cases. The investment costs and $W_{t}$ values are taken from \cite{devine2024analysing}, while all other parameters' values are taken from \cite{devine2023cournot}. The parameter values can be found in Online Appendix \footnote{https://figshare.com/articles/dataset/Initial-Data-Irish-El-Data/28785299}.
\vspace{-1em}

\section{Results and Conclusions}
\label{sec:5}
\vspace{-1em}
Table \ref{tab1:electricity_market} summarizes the results of the three proposed models for different demand cases in Ireland for 2026. The results of this study highlight differences in electricity generation, CO$_2$ emissions, investment levels, and market prices under different market structures in Ireland’s electricity market. 
For each test case, the output demand under perfect competition corresponds with the demand levels projected in Ireland's Eirgrid report, whereas the imperfect competition model suggests a slightly lower generation level. This difference is attributed to strategic behavior by firms in the imperfect market setting; that is, strategic players induce a higher price and hence lower demand.

\begin{table}
    \caption{Key results of proposed models for different electricity demand cases in Ireland 2026}
    \label{tab1:electricity_market}
    \centering
    \resizebox{10cm}{!}{%
    \begin{tabular}{lccc}
        \toprule
        \textbf{Demand cases} & \textbf{Low demand} & \textbf{Median demand} & \textbf{High demand} \\
        \midrule
        \multicolumn{4}{l}{\textbf{Total electricity generation (TWh)}} \\
        Ireland (Eirgrid) test cases & 44.70 & 48.80 & 52.40 \\
        Perfect competition generation & 44.71 & 48.93 & 52.68 \\
        Perfect with UC constraint generation & 43.83 & 48.15 & 51.97 \\
        Imperfect competition generation & 39.72 & 43.36 & 46.61 \\
        \midrule
        \multicolumn{4}{l}{\textbf{Total CO\textsubscript{2} emission (MTons)}} \\
        Perfect competition & 10.37 & 11.29 & 12.02 \\
        Perfect competition with UC constraint & 10.90 & 11.63 & 12.27 \\
        Imperfect competition & 8.27 & 9.11 & 9.81 \\
        \midrule
        \multicolumn{4}{l}{\textbf{CO\textsubscript{2} emission per TWh (MTons/TWh)}} \\
        Perfect competition & 4.31 & 4.33 & 4.38 \\
        Perfect competition with UC constraint & 4.02 & 4.14 & 4.24 \\
        Imperfect competition & 4.80 & 4.76 & 4.75 \\
        \midrule
        \multicolumn{4}{l}{\textbf{Total investment (MW)}} \\
        Perfect competition & 0 & 577 & 1240 \\
        Perfect competition with UC constraint & 3486 & 4789 & 5901 \\
        Imperfect competition & 167 & 525 & 1026 \\
        \midrule
        \multicolumn{4}{l}{\textbf{Contribution of renewables technologies (\%)}} \\
        Perfect competition & 54.86 & 54.41 & 54.53 \\
        Perfect competition with UC constraint & 59.61 & 59.36 & 59.14 \\
        Imperfect competition & 59.24 & 58.16 & 57.78 \\
        \midrule
        \multicolumn{4}{l}{\textbf{Average market price of electricity (€ / MWh)}} \\
        Perfect competition & 49.49 & 52.10 & 53.70 \\
        Perfect competition with UC constraint & 63.87 & 64.70 & 65.07 \\
        Imperfect competition & 127.14 & 138.80 & 148.32 \\
        \bottomrule
    \end{tabular}
    }
\end{table}

Additionally, the CO$_2$ emissions per TWh vary across models, with the imperfect competition model exhibiting lower overall emissions despite lower total generation. This suggests that under an imperfect market structure, firms may prioritize more efficient or lower-carbon technologies in response to market constraints. Conversely, the perfect competition model, particularly with unit commitment constraints, shows a higher renewable contribution, indicating that in a competitive setting, investment in renewable technologies might be more favorable. However, this scenario also results in higher electricity prices, reflecting the cost implications of incorporating unit-level operational constraints.

These findings underscore the trade-offs inherent in different market structures and modeling choices. While perfect competition encourages increased renewable penetration, it may not fully capture the CO$_2$ emissions reduction that an imperfect competition model can reflect. Similarly, investment levels in generation capacity are higher in the perfect competition framework with unit commitment constraints, with investments reaching up to 5901 MW in the high-demand test case. The imperfect competition model, in contrast, suggests a more measured investment trajectory, with a maximum investment of 1026 MW. These variations imply that policy decisions regarding market design should carefully balance competition, emissions targets, and investment sustainability to achieve Ireland’s long-term energy objectives. 

In addition, the market price of electricity under imperfect competition is significantly higher (up to 148.32 €/MWh in high demand) compared to perfect competition scenarios (53.70 €/MWh without unit commitment constraints), highlighting the cost implications of different market structures. Moreover, choosing different modeling assumptions will lead to different projections.

In conclusion, this study provides insights for Ireland’s electricity market planning by analyzing both perfect and imperfect competition structures. The key and novel contribution is the integration of long-term investment decisions with detailed Irish market data, along with a comparative study of a unit commitment model. This approach allows for a more realistic assessment of market dynamics, investment behavior, CO$_2$ emissions, and electricity prices. The results highlight potential outcomes under different market conditions. Future research could further refine these models by incorporating additional factors such as policy incentives and technological advancements to enhance the effectiveness of investment planning in Ireland’s evolving electricity market.
\vspace{-2em}
\section*{Acknowledgments}
\vspace{-1em}
\footnotesize
D. Hosseinnezhad and M. T. Devine acknowledge funding from Research Ireland and co-funding partners under grant number 21/SPP/3756 through the NexSys Strategic Partnership Programme. 
\vspace{-1em}
%
%
%

\vspace{-0.5em}
\end{document}